\documentclass[12pt]{article}
\usepackage{graphicx}

\usepackage{epsfig}
\usepackage{amsmath}
\usepackage{amssymb}
\usepackage{amsthm}
\usepackage{latexsym}
\usepackage{cite}
\usepackage{amsbsy,amsfonts,euscript,eepic,rotating}

\newcommand{\comment}[1]{}
\newtheorem{theorem}{Theorem}
\newtheorem{lemma}{Lemma}[section]

\newtheorem{hypothesis}{Hypothesis}[section]
\begin{document}

\title{\LARGE
{\bf Scaling Limits of Two-Dimensional Percolation: an Overview}
}

\author{
{\bf Federico Camia}
\thanks{Research supported in part by a Veni grant of the
NWO (Dutch Organization for Scientific Research).}\,
\thanks{E-mail: fede@few.vu.nl}\\
{\small \sl Department of Mathematics, Vrije Universiteit Amsterdam}\\
{\small \sl De Boelelaan 1081a, 1081 HV Amsterdam, The Netherlands}
}

\date{}

\maketitle

\begin{abstract}
We present a review of the recent progress on percolation scaling limits
in two dimensions. In particular, we will consider the convergence of
critical crossing probabilities to Cardy's formula and of the critical
exploration path to chordal SLE$_6$, the full scaling limit of critical
cluster boundaries, and near-critical scaling limits.
\end{abstract}



\section{Introduction}

Percolation as a mathematical theory was introduced by Broadbent and Hammersley~\cite{broadbent,bh}
to model the spread of a gas or a fluid through a porous medium.
To mimic the randomness of the medium, they declared the edges of the $d$-dimensional cubic lattice
independently \emph{open} (to the passage of the gas or fluid) with probability $p$ or \emph{closed}
with probability $1-p$. Since then, many variants of this simple model have been studied, attracting
the interest of both mathematicians and physicists.

Mathematicians are interested in percolation because of its deceiving simplicity which hides difficult
and elegant results.
From the point of view of physicists, percolation is maybe the simplest statistical mechanical model
undergoing a continuous phase transition as the value of the parameter $p$ is varied, with all the
standard features typical of critical phenomena (scaling laws, conformal invariance, universality).
On the applied side, percolation has been used to model the spread of a disease, a fire or a rumor,
the displacement of oil by water, the behavior of random electrical circuits, and more recently the
connectivity properties of communication networks.

We will concentrate on a version of the model in which each vertex of the triangular
lattice, identified by duality with the corresponding face of the hexagonal lattice
(see Figure~\ref{fig-tri-hex}), is independently colored black with probability $p$
or white otherwise.
Questions regarding the geometry of this random coloring (for instance, whether
there exists a path on black sites connecting the opposite edges of a given rectangle)
can be expressed in terms of the behavior of {\em clusters} (i.e., maximal connected
monochromatic subsets of the lattice) or of the boundaries between them (we will
sometimes call such boundaries or portions of them {\em percolation interfaces}).

It is well known (see, e.g., \cite{kesten-book,grimmett-book,br-book}) that in this
model if $p>1/2$ (respectively, $p<1/2$) there is an infinite black (resp., white)
cluster, while for $p=1/2$ there is no infinite cluster of either color. The latter
value is the critical threshold $p_c$ of the model, at which the percolation phase
transition occurs. As testified by~\cite{kesten-book,grimmett-book}, we have had for
some time a good understanding of the \emph{subcritical} ($p<p_c$) and \emph{supercritical}
($p>p_c$) phases.
As for the \emph{critical behavior} ($p$ equal to or approaching $p_c$), despite
some important achievements (see, in particular, \cite{kesten} and~\cite{kesten-book,grimmett-book}
as general references), a complete and rigorous understanding seemed out of reach for
any two-dimensional percolation model until the introduction of the Stochastic Loewner
Evolution (SLE) by Oded Schramm~\cite{schramm} and the proof of conformal invariance
by Stanislav Smirnov~\cite{smirnov1,smirnov1-long}.

The percolation phase transition is a purely geometric transition that physicists have
successfully studied with the methods of continuous phase transitions, or critical
phenomena. In the theory of critical phenomena it is usually assumed that a physical
system near a continuous phase transition is characterized by a single length scale
(the \emph{correlation length}) in terms of which all other lengths should be measured.
When combined with the experimental observation that the correlation length diverges at
the phase transition, this simple but strong assumption, known as the scaling hypothesis,
leads to the belief that at criticality the system has no characteristic length, and is
therefore invariant under scale transformations. This implies that all thermodynamic
functions at criticality are homogeneous functions, and predicts the appearance of
power laws.

It also suggests that for models of critical systems realized on a lattice,
one can attempt to take a \emph{continuum scaling limit} in which the mesh
of the lattice is sent to zero while focus is kept on ``macroscopic"
observables that capture the large scale behavior.
In the limit, the discrete model should converge to a continuum one that
encodes the large scale properties of the original model, containing at
the same time more symmetry. In many cases, this allows to derive extra
insight by combining methods of discrete mathematics with considerations
inspired by the continuum limit picture. The simplest
example of such a continuum random model is Brownian motion,
which is the scaling limit of the simple random walk.
In general, though, the complexity of the discrete model makes
it impossible to even guess the nature of the scaling limit,
unless some additional feature can be shown to hold, which can
be used to pin down properties of the continuum limit.
Two-dimensional critical systems belong to the class of models
for which this can be done, and the additional feature is
\emph{conformal invariance}, as predicted by physicists since the
early seventies~\cite{polyakov1,polyakov2}.

The connection between the scaling limit of critical percolation interfaces
(i.e., boundaries between clusters of different colors) and SLE has
led to tremendous progress in recent years, not only providing a rigorous
derivation of many of the results previously obtained by physicists, but also
deepening our geometric understanding of critical percolation, and critical
phenomena in general. The main power of SLE stems from the fact that it allows
to compute different quantities; for example, percolation crossing probabilities
and various percolation critical exponents. In general, relating the scaling
limit of a critical lattice model to SLE allows for a rigorous determination of
some aspects of the large scale behavior of the lattice model. It also provides
deeper insight into geometric aspects that are not easily accessible with the
methods developed by physicists to study critical phenomena. For mathematicians,
the biggest advantage of SLE over those methods lies maybe in its mathematical
rigor. However, many physicists working on critical phenomena have promptly
recognized the importance of SLE and added it to their toolbox.

\section{SLE and CLE} \label{sec-sle-cle}

The Stochastic Loewner Evolution (also called Schramm-Loewner Evolution)
with parameter $\kappa > 0$ ($\text{SLE}_{\kappa}$) was introduced by
Schramm~\cite{schramm} as a tool for studying the large scale behavior
of two-dimensional discrete (defined on a lattice) probabilistic models
whose scaling limits are expected to be conformally invariant.
In this section we define the chordal version of $\text{SLE}_{\kappa}$;
for more on the subject, the interested reader can consult
the original paper~\cite{schramm} as well as~\cite{rs}, the fine
reviews by Lawler~\cite{lawler1}, Kager and Nienhuis~\cite{kn},
and Werner~\cite{werner4}, and Lawler's book~\cite{lawler2}.

Let $\mathbb H$ denote the upper half-plane.
For a given continuous real function $U_t$ with $U_0 = 0$,
define, for each $z \in \overline{\mathbb H}$, the function
$g_t(z)$ as the solution to the ODE
\begin{equation}
\partial_t g_t(z) = \frac{2}{g_t(z) - U_t},
\end{equation}
with $g_0(z) = z$.
This is well defined as long as $g_t(z) - U_t \neq 0$,
i.e., for all $t < T(z)$, where
\begin{equation}
T(z) := \sup \{ t \geq 0 : \min_{s \in [0,t]} | g_s(z) - U_s| > 0 \}.
\end{equation}
Let $K_t := \{ z \in \overline{\mathbb H} : T(z) \leq t \}$
and ${\mathbb H}_t := {\mathbb H} \setminus K_t$; it can be shown that
$K_t$ is bounded, ${\mathbb H}_t$ is simply connected, and $g_t$ is a
conformal map from ${\mathbb H}_t$ onto $\mathbb H$.
For each $t$, it is possible to write $g_t(z)$ as
\begin{equation}
g_t(z) = z + \frac{2t}{z} + O\left(\frac{1}{z^2}\right),
\end{equation}
when $z \to \infty$.
The family $(K_t, t \geq 0)$ is called the {\em Loewner chain}
associated to the driving function $(U_t, t \geq 0)$.

We will call {\em chordal $\text{SLE}_{\kappa}$} the Loewner chain $(K_t, t \geq 0)$
that is obtained when the driving function $U_t = \sqrt{\kappa} B_t$ is $\sqrt{\kappa}$
times a standard real-valued Brownian motion $(B_t, t \geq 0)$ with $B_0 = 0$.
(Note that $\text{SLE}_{\kappa}$ is often defined as the family of conformal maps
$(g_t, t \geq 0)$, but we find the above definition more convenient for our purposes.)

For all $\kappa \geq 0$, chordal $\text{SLE}_{\kappa}$ is almost surely generated
by a continuous random curve $\gamma$ in the sense that, for all $t \geq 0$,
${\mathbb H}_t = {\mathbb H} \setminus K_t$ is the unbounded connected
component of ${\mathbb H} \setminus \gamma[0,t]$; $\gamma$ is called the
{\em trace} of chordal $\text{SLE}_{\kappa}$.

It is not hard to see, as first argued by Schramm~\cite{schramm}, that a
continuous random curve $\gamma$ in the upper half-plane starting at the origin
and going to infinity must be an SLE trace if it possesses the following
{\em conformal Markov property} (sometimes called \emph{domain Markov property}).
For any fixed $T \in {\mathbb R}$, conditioned on $\gamma[0,T]$, the image
under $g_T$ of $\gamma[T,\infty)$ is distributed like an independent copy
of $\gamma$, up to a time reparametrization.
This implies that the driving function $U_t$ in the Loewner chain associated to
the curve $\gamma$ is continuous and has stationary and independent increments.
If the time parametrization implicit in the definition of chordal SLE$_k$ and the
discussion preceding it is chosen for $\gamma$, then scale invariance also implies
that the law of $U_t$ is the same as the law of $\lambda^{-1/2} U_{\lambda t}$ when
$\lambda>0$. These properties together imply that $U_t$ must be a constant multiple
of standard Brownian motion.

Let now $D \subset {\mathbb C}$ ($D \neq {\mathbb C}$) be a simply
connected domain whose boundary is a continuous curve.
By Riemann's mapping theorem, there are (infinitely many) conformal
maps from the upper half-plane $\mathbb H$ onto $D$.
In particular, given two distinct points $a,b \in \partial D$
(or more accurately, two distinct prime ends), there exists a
conformal map $f$ from $\mathbb H$ onto $D$ such that $f(0)=a$
and $f(\infty) \equiv \lim_{|z| \to \infty} f(z) = b$.
In fact, the choice of the points $a$ and $b$ on the boundary
of $D$ only characterizes $f(\cdot)$ up to a scaling factor $\lambda>0$,
since $f(\lambda \cdot)$ would also do.

Suppose that $(K_t, t \geq 0)$ denotes chordal $\text{SLE}_{\kappa}$ in
$\mathbb H$ as defined above; we define chordal $\text{SLE}_{\kappa}$
$(\tilde K_t, t \geq 0)$ in $D$ from $a$ to $b$ as the image of the
Loewner chain $(K_t, t \geq 0)$ under $f$.
It is possible to show, using scaling properties of
$\text{SLE}_{\kappa}$, that the law of $(\tilde K_t, t \geq 0)$
is unchanged, up to a linear time-change, if we replace
$f(\cdot)$ by $f(\lambda \cdot)$.
This makes it natural to consider $(\tilde K_t, t \geq 0)$ as
a process from $a$ to $b$ in $D$, ignoring the role of $f$.
The trace of chordal SLE in $D$ from $a$
to $b$ will be denoted by $\gamma_{D,a,b}$.


We now move from the conformally invariant random curves of SLE to
collections of conformally invariant random loops, and introduce the concept
of Conformal Loop Ensemble (CLE --- see~\cite{werner3,werner5,sheffield}).
The key feature of a CLE is a sort of ``{\em conformal restriction/renewal property}."
Roughly speaking, a CLE in $D$ is a random collection ${\cal L}_D$ of loops
such that if all the loops intersecting a (closed) subset of $D$ or of
its boundary are removed, the loops in any one of the various remaining
(disjoint) subdomains of $D$ form a random collection of loops distributed
as an independent copy of ${\cal L}_D$ conformally mapped to that subdomain
(see Theorem~\ref{thm-cle}).
We will not attempt to be more precise here since somewhat different definitions
(although, in the end, substantially equivalent) have appeared in the literature,
but the meaning of the conformal restriction/renewal property should be clear
from Theorem~\ref{thm-cle}.

For formal definitions and more discussion on the properties of a CLE,
the reader is referred to the original literature on the subject~\cite{werner3,werner5,sheffield},
where it is shown that there is a one-parameter family $\text{CLE}_{\kappa}$
of conformal loop ensembles with the above conformal restriction/renewal
property and that for $\kappa \in (8/3,8]$, the $\text{CLE}_{\kappa}$ loops
locally ``look like" $\text{SLE}_{\kappa}$ curves.

There are numerous lattice models that can be described in terms of random
curves and whose scaling limits are assumed (and in a few cases proved) to
be conformally invariant.
These include the Loop Erased Random Walk, the Self-Avoiding Walk and the
Harmonic Explorer, all of which can be defined as polygonal paths along
the edges of a lattice.
The Ising, Potts and percolation models instead are naturally defined in
terms of clusters, and the interfaces between different clusters form random
loops.
In the $O(n)$ model, configurations of loops along the edges of the hexagonal
lattice are weighted according to the total number and length of the loops.
All of these models are supposed to have scaling limits described by
$\text{SLE}_{\kappa}$ or $\text{CLE}_{\kappa}$ for some value of $\kappa$
between $2$ and $8$.
For more information on these lattice models and their scaling limits, the
interested reader can consult~\cite{kn,cardy2,cardy3,werner5,sheffield}.
In the case of percolation, the connection with $\text{SLE}_6$ and $\text{CLE}_6$
has been made rigorous~\cite{smirnov1,smirnov1-long,cn2,cn3,cn4}.

\section{Definitions and Preliminary Results} \label{sec-prelim}

In the rest of the paper we will consider critical site percolation on the
triangular lattice, for which conformal invariance in the scaling limit has
been rigorously proved~\cite{smirnov1,smirnov1-long}.
A precise formulation of conformal invariance, attributed to Michael Aizenman,
is that the probability that a percolation cluster crosses between two disjoint
segments of the boundary of some simply connected domain should converge to a
conformally invariant function of the domain and the two segments of the boundary.
This conjecture is connected with the extensive numerical investigations
reported in~\cite{lps}.
A formula for the purposed limit was then derived~\cite{cardy1} by John Cardy
using (non-rigorous) 
computations that assumed the existence of a conformally invariant scaling limit.
The interest of mathematicians was already evident in~\cite{lps}, but a proof
of the conjecture~\cite{smirnov1,smirnov1-long} (and of Cardy's formula) did
not come until 2001.

\begin{figure}[!ht]
\begin{center}
\includegraphics[width=5cm]{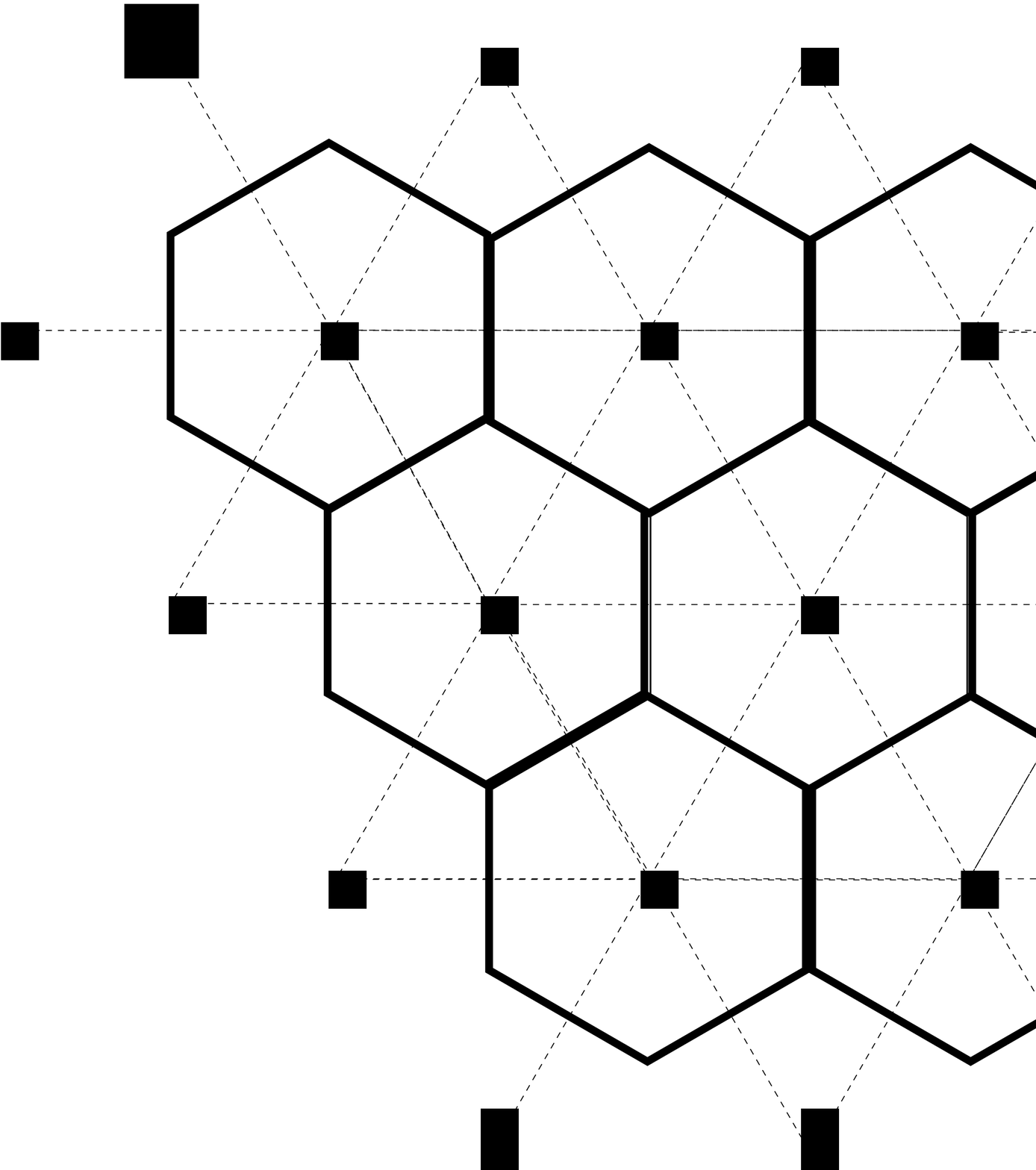}
\caption{Embedding of the triangular and hexagonal lattices in ${\mathbb R}^2$.}
\label{fig-tri-hex}
\end{center}
\end{figure}

We will denote by $\cal T$ the two-dimensional triangular lattice,
whose sites are identified with the elementary cells of a regular hexagonal
lattice $\cal H$ embedded in the plane as in Fig.~\ref{fig-tri-hex}.
We say that two hexagons are neighbors (or that they are adjacent) if
they have a common edge.
A sequence $(\xi_0, \ldots, \xi_n)$ of 
hexagons of $\cal H$ such that $\xi_{i-1}$ and $\xi_i$
are neighbors for all $i= 1, \ldots, n$ and $\xi_i \neq \xi_j$
whenever $i \neq j$ will be called a {\em $\cal T$-path}.
A $\cal T$-path whose first and last hexagons are neighbors will be
called a {\em $\cal T$-loop}.

\boldmath
\subsection{Compactification of ${\mathbb R}^2$}
\unboldmath

When taking the scaling limit, as the lattice spacing
$\delta \to 0$, one can focus on fixed finite regions,
$\Lambda \subset {\mathbb R}^2$, or consider the whole
${\mathbb R}^2$ at once.
The second option avoids dealing with boundary conditions,
but requires an appropriate choice of metric.

A convenient way of dealing with the whole ${\mathbb R}^2$
is to replace the Euclidean metric with a distance function
$\Delta(\cdot,\cdot)$ defined on
${\mathbb R}^2 \times {\mathbb R}^2$ by
\begin{equation} \nonumber
\Delta(u,v) := \inf_{\varphi} \int (1 + | {\varphi} |^2)^{-1} \, ds,
\end{equation}
where the infimum is over all smooth curves $\varphi(s)$
joining $u$ with $v$, parameterized by arclength $s$, and
where $|\cdot|$ denotes the Euclidean norm.
This metric is equivalent to the Euclidean metric in bounded
regions, but it has the advantage of making ${\mathbb R}^2$
precompact.
Adding a single point at infinity yields the compact space
$\dot{\mathbb R}^2$ which is isometric, via stereographic
projection, to the two-dimensional sphere.

\subsection{The space of curves} \label{space}

In dealing with the scaling limit we use the approach of
Aizenman-Burchard~\cite{ab}.
Let $D$ be a Jordan domain and denote by ${\cal S}_D$ the
complete separable metric space of continuous curves in the
closure $\overline D$ of $D$ with the metric~(\ref{distance})
defined below.
Curves are regarded as equivalence classes of continuous
functions from the unit interval to $\overline D$, modulo
monotonic reparametrizations.
$\gamma$ will represent a particular curve and $\gamma(t)$ a
parametrization of $\gamma$; ${\cal F}$ will represent a set
of curves (more precisely, a closed subset of ${\cal S}_D$).
$\text{d}(\cdot,\cdot)$ will denote the uniform metric
on curves, defined by
\begin{equation} \label{distance}
\text{d} (\gamma_1,\gamma_2) := \inf
\sup_{t \in [0,1]} |\gamma_1(t) - \gamma_2(t)|,
\end{equation}
where the infimum is over all choices of parametrizations
of $\gamma_1$ and $\gamma_2$ from the interval $[0,1]$.
The distance between two closed sets of curves is defined
by the induced Hausdorff metric as follows:
\begin{equation} \label{hausdorff}
\text{dist}({\cal F},{\cal F}') \leq \varepsilon
\Leftrightarrow (\forall \, \gamma \in {\cal F}, \, \exists \,
\gamma' \in {\cal F}' \text{ with }
\text{d} (\gamma,\gamma') \leq \varepsilon,
\text{ and vice versa}).
\end{equation}
The space $\Omega_D$ of closed subsets of ${\cal S}_D$
(i.e., collections of curves in $\overline D$) with the
metric~(\ref{hausdorff}) is also a complete separable metric
space. 
For each fixed $\delta>0$, the random curves that we consider are
polygonal paths on the edges of the hexagonal lattice $\delta{\cal H}$,
dual to the triangular lattice $\delta{\cal T}$, corresponding to
boundaries between black and white clusters.

We will also consider the complete separable metric space ${\cal S}$
of continuous curves in $\dot{\mathbb R}^2$ with the distance
\begin{equation} \label{Distance}
\text{D} (\gamma_1,\gamma_2) := \inf
\sup_{t \in [0,1]} \Delta(\gamma_1(t),\gamma_2(t)),
\end{equation}
where the infimum is again over all choices of parametrizations
of $\gamma_1$ and $\gamma_2$ from the interval $[0,1]$.
The distance between two closed sets of curves is again
defined by the induced Hausdorff metric as follows:
\begin{equation} \label{hausdorff-D}
\text{Dist}({\cal F},{\cal F}') \leq \varepsilon
\Leftrightarrow (\forall \, \gamma \in {\cal F}, \, \exists \,
\gamma' \in {\cal F}' \text{ with }
\text{D} (\gamma,\gamma') \leq \varepsilon
\text{ and vice versa}).
\end{equation}
The space $\Omega$ of closed sets of $\cal S$
(i.e., collections of curves in $\dot{\mathbb R}^2$)
with the metric~(\ref{hausdorff-D}) is also a complete
separable metric space.

When we talk about convergence in distribution of random curves,
we always mean with respect to the uniform metrics~(\ref{distance})
or~(\ref{Distance}), while when we deal with closed collections
of curves, we always refer to the metrics~(\ref{hausdorff})
or~(\ref{hausdorff-D}).
In this paper, the space $\Omega$ of closed sets of $\cal S$ is
used for collections of boundary contours and their scaling limits.

\subsection{Existence of subsequential scaling limits}

Aizenman and Burchard \cite{ab} provided a crucial step in the analysis
of continuum scaling limits for a very general class of models, including
percolation. They introduced the metric space of curves described earlier,
and formulated a condition which is sufficient for the existence of
subsequential scaling limits (via a compactness argument). In order to state
the condition, we need some definitions. For $\delta >0$, let $\mu_{\delta}$
be any probability measure supported on collections of curves that are
polygonal paths on the edges of the scaled hexagonal lattice $\delta{\cal H}$.
For $x \in {\mathbb R}^2$ and $R>r>0$, let $A(x;r,R)$ be the annulus with
inner radius $r$ and outer radius $R$ centered at $x$. We say that an annulus
is crossed by a curve if the curve intersects both the inner and the outer
circles.

\begin{hypothesis}\label{Hypothesis 1}
For all $k < \infty$ and for all annuli $A(x;r,R)$ with $\delta \leq r \leq R\leq 1$,
the following bound holds uniformly in $\delta$:

\begin{equation} \nonumber
\mu_{\delta}\left(A(x;r,R) \mbox{ is crossed by } k \mbox{ disjoint curves}\right) \leq K_k
\left(\frac{r}{R}\right)^{\phi(k)}
\end{equation}
for some $K_k < \infty$ and $\phi(k) \rightarrow \infty$ as $k \rightarrow \infty$.
\end{hypothesis}

If the above hypothesis is satisfied, it follows from~\cite{ab} that for every
sequence $\delta_j \downarrow 0$, there exists a subsequence $\{\delta_{j_i}\}_{i \in {\mathbb N}}$
such that $\{\mu_{\delta_{j_i}}\}_{i \in {\mathbb N}}$ has a limit. More precisely, we have the
following theorem, which follows from a more general result proved in~\cite{ab}.

\begin{theorem} \label{Thm Aizenman}
Hypothesis \ref{Hypothesis 1} implies that for any sequence $\delta_j \downarrow 0$,
there exists a subsequence $\{\delta_{j_{i}}\}_{i \in {\mathbb N}}$ and a probability
measure $\mu$ on $\Omega$ such that $\mu_{\delta_{j_{i}}}$ converges weakly to $\mu$
as $i \rightarrow \infty$.
\end{theorem}

It was already remarked in~\cite{ab} that the above hypothesis can be verified,
using the RSW theorem~\cite{russo,sewe} and the BK inequality~\cite{vdbk}, for
two-dimensional critical and near-critical (see Section~\ref{sec-near-crit})
percolation. The same conclusion follows from results in~\cite{nolin},
and is obtained in Proposition~1 of~\cite{nowe}. In fact, the more general
result below can be easily proved (see~\cite{cjm}). Such generality is useful
when considering near-critical percolation (see Sect.~\ref{sec-near-crit}
and~\cite{cjm}). Note that the sequence $\{p_j\}_{j \in {\mathbb N}}$ in
the lemma below need not converge, but the lemma is typically applied to
converging sequences.

\begin{lemma} \label{lemma-hypo}
Let $\{\mu_{\delta_j,p_j}\}_{j \in {\mathbb N}}$ be a sequence of measures on
boundary contours induced by percolation on $\delta_j{\cal T}$ with parameters $p_j$.
For {\em any} sequence $\delta_j \to 0$ and {\em any} choice of the collection
$\{p_j\}_{j \in {\mathbb N}}$, Hypothesis~\ref{Hypothesis 1} holds.
\end{lemma}

Theorem~\ref{Thm Aizenman} and Lemma~\ref{lemma-hypo} guarantee the existence
of subsequential scaling limits of the measure $\mu_{\delta,p}$ on boundary
contours induced by percolation on $\delta{\cal T}$ with parameter $p$ which
can depend on $\delta$. But what about uniqueness? And what properties of
scaling limits can be proved? Exponential decay of cluster sizes away from
$p_c$ implies that the scaling limit of boundary contours is trivial (i.e.,
consists of no curves with diameter larger than zero) if $p$ does not tend
to $p_c$. Therefore, the interesting cases are $p=p_c$ for all $\delta$,
and $p \to p_c$ as $\delta \to 0$. In the next section we will consider
the first case.

\section{Conformal Invariance of Critical Percolation} \label{sec-conformal}

In this section we consider the percolation model on the triangular lattice introduced
earlier, and let $p=p_c=1/2$. The special feature of this case that will allow us to
determine the scaling limit is {\em conformal invariance}.

\subsection{Cardy's Formula and Smirnov's Theorem}

Let $D \in {\mathbb R}^2$ be a bounded, simply connected domain containing
the origin whose boundary $\partial D$ is a continuous curve. We will mostly
be concerned with {\em Jordan domains}, i.e., bounded domains whose boundaries
are simple, closed curves.
Let $\phi:\overline{\mathbb D} \to \overline D$ be the (unique) continuous
function that maps the unit disc $\mathbb D$ onto $D$ conformally and such
that $\phi(0)=0$ and $\phi'(0)>0$.
Let $z_1,z_2,z_3,z_4$ be four points of $\partial D$ in counterclockwise
order --- i.e., such that $z_j=\phi(w_j), \,\,\, j=1,2,3,4$, with
$w_1,\ldots,w_4$ in counterclockwise order.
Also, let $\eta = \frac{(w_1-w_2)(w_3-w_4)}{(w_1-w_3)(w_2-w_4)}$.
Cardy's formula~\cite{cardy1} for the probability $\Phi_{D}(z_1,z_2;z_3,z_4)$
of a ``crossing" inside $D$ from the counterclockwise arc $\overline{z_1 z_2}$
to the counterclockwise arc $\overline{z_3 z_4}$ is
\begin{equation} \label{cardy-formula}
\Phi_{D}(z_1,z_2;z_3,z_4) =
\frac{\Gamma(2/3)}{\Gamma(4/3) \Gamma(1/3)} \eta^{1/3} {}_2F_1(1/3,2/3;4/3;\eta),
\end{equation}
where ${}_2F_1$ is a hypergeometric function.

For a given mesh $\delta>0$, the probability of a black crossing inside $D$
from the counterclockwise arc $\overline{z_1 z_2}$ to the counterclockwise
arc $\overline{z_3 z_4}$ is the probability of the existence of a
$\cal T$-path $(\xi_0,\ldots,\xi_n)$ such that $\xi_0$ intersects the
counterclockwise arc $\overline{z_1 z_2}$, $\xi_n$ intersects the
counterclockwise arc $\overline{z_3 z_4}$, $\xi_1,\ldots,\xi_{n-1}$
are all contained in $D$, and $\xi_0,\ldots,\xi_n$ are all black.
Smirnov~\cite{smirnov1,smirnov1-long} proved that crossing probabilities
converge in the scaling limit to conformally invariant functions of the
domain and the four points on its boundary, and identified the limit
with Cardy's formula~(\ref{cardy-formula}).

One formulation of the result, somewhat less general than
Smirnov's~\cite{smirnov1,smirnov1-long} in terms of the domains considered,
is given below. A very detailed (and rather lengthy) proof of the theorem
can be found in~\cite{br-book} (see also~\cite{beffara}).
\begin{theorem} 
\label{cardy-smirnov}
Let $D$ be a Jordan domain of the plane.
As $\delta \to 0$, the limit of the probability of a black crossing inside $D$
from the counterclockwise arc $\overline{z_1 z_2}$ to the counterclockwise arc
$\overline{z_3 z_4}$ exists, is a conformal invariant of $(D,z_1,z_2,z_3,z_4)$,
and is given by Cardy's formula~(\ref{cardy-formula}).
\end{theorem}

The proof of Smirnov's theorem is based on the identification of certain
generalized crossing probabilities that are almost discrete harmonic functions
and whose scaling limits converge to harmonic functions.
The behavior on the boundary of such functions is easy to determine and is
sufficient to specify them uniquely.
The relevant crossing probabilities can be expressed in terms of the boundary
values of such harmonic functions and, as a consequence, are invariant under
conformal transformations of the domain and the two segments of its boundary.

The presence of a black crossing in $D$ from the counterclockwise boundary arc
$\overline{z_1 z_2}$ to the counterclockwise boundary arc $\overline{z_3 z_4}$
can be determined using a clever algorithm that explores the percolation
configuration inside $D$ starting at, say, $z_1$ and assuming that the
hexagons just outside $D$ along $\overline{z_1 z_2}$ are all black and
those along $\overline{z_4 z_1}$ are all white (a more precise definition
of the algorithm can be formulated in the context of the {\em lattice domains}
introduced in the next section).
The exploration proceeds following the interface between the black
cluster adjacent to $\overline{z_1 z_2}$ and the white cluster adjacent
to $\overline{z_4 z_1}$.
A black crossing is present if the exploration process reaches $\overline{z_3 z_4}$
before $\overline{z_2 z_3}$.
This {\em exploration process} and the {\em exploration path}
associated to it were introduced by Schramm in~\cite{schramm}.

The connection between crossing probabilities and the percolation exploration
process suggests that maybe the exploration path should also be conformally
invariant in some appropriate sense. A precise conjecture was formulate by
Schramm~\cite{schramm}, and will be discussed in the next section.


\subsection{Convergence of the Percolation Exploration Path to SLE$_6$}

The exploration process can be carried out in ${\mathbb H} \cap {\cal H}$,
where the hexagons in the lowest row and to the left of a chosen hexagon
have been colored yellow and the remaining hexagons in the lowest row have
been colored blue.
This produces an infinite exploration path, whose scaling limit was
conjectured~\cite{schramm} by Schramm to converge to $\text{SLE}_6$.

It is easy to see that the exploration process is Markovian in the sense
that, conditioned on the exploration up to a certain (stopping) time,
the future of the exploration evolves in the same way as the past except
that it is now performed in a different domain. The new domain is the
upper half-plane minus the explored region, and some of the explored
hexagons are now part of the new boundary. 
This observation, combined with the connection between the exploration
process and crossing probabilities, Smirnov's theorem about the conformal
invariance of crossing probabilities in the scaling limit, and Schramm's
characterization of SLE via the conformal Markov property discussed in
Sect.~\ref{sec-sle-cle}, strongly supports the above conjecture.


As we now explain, the natural setting to define the exploration process
is that of {\em lattice domains}, i.e., sets $D^{\delta}$ of hexagons of
$\delta\cal H$ that are {\em connected} in the sense that any two hexagons
in $D^{\delta}$ can be joined by a $(\delta\cal T)$-path contained in $D^{\delta}$.
We say that a bounded lattice domain $D^{\delta}$ is {\em simply connected}
if both $D^{\delta}$ and $\delta{\cal T} \setminus D^{\delta}$ are connected.
We call $D^{\delta}$ a {\em lattice-Jordan} domain if it is a bounded,
simply connected lattice domain such that the hexagons adjacent to
$D^{\delta}$ form a $(\delta\cal T)$-loop.

Given a lattice-Jordan domain $D^{\delta}$, the set of hexagons
adjacent to $D^{\delta}$ can be partitioned into two (lattice-)connected
sets. If those two sets of hexagons are assigned different colors,
for any coloring of the hexagons inside $D^{\delta}$, there is an
interface between two clusters of different colors starting and
ending at two boundary points, $a^{\delta}$ and $b^{\delta}$,
corresponding to the locations on the boundary of $D^{\delta}$
where the color changes.
If one performs an exploration process in $D^{\delta}$ starting
at $a^{\delta}$, one ends at $b^{\delta}$, producing an exploration
path $\gamma^{\delta}$ that traces the entire interface from
$a^{\delta}$ to $b^{\delta}$.

Given a planar domain $D$, we denote by $\partial D$ its topological
boundary. Assume that $D$ is a Jordan domain that contains the origin.
Given $a,b \in \partial D$, we write $(D^{\delta},a^{\delta},b^{\delta}) \to (D,a,b)$
if $\partial D^{\delta}$ converges to $\partial D$ in the metric~(\ref{distance})
and $a^{\delta},b^{\delta} \in \partial D^{\delta}$ converge respectively
to $a,b \in \partial D$ as $\delta \to 0$. A detailed proof of the next theorem,
first stated in a slightly different form in~\cite{smirnov1,smirnov1-long},
can be found in~\cite{cn3}.


\begin{theorem} \label{thm-conv-to-sle}
Let $(D,a,b)$ be a Jordan domain with two distinct selected points
on its boundary $\partial D$.
Then, for lattice-Jordan domains $D^{\delta}$ from $\delta{\cal H}$ with
$a^{\delta},b^{\delta} \in \partial D^{\delta}$ such that
$(D^{\delta},a^{\delta},b^{\delta}) \to (D,a,b)$ as $\delta \to 0$,
the percolation exploration path $\gamma^{\delta}_{D,a,b}$ in $D^{\delta}$
from $a^{\delta}$ to $b^{\delta}$ converges in distribution to the
trace $\gamma_{D,a,b}$ of chordal $\text{\emph{SLE}}_6$ in $D$
from $a$ to $b$, as $\delta \to 0$.
\end{theorem}

\subsection{Critical Exponents} \label{sec-crit-exp}

The convergence of the percolation exploration path to SLE$_6$, combined with
Kesten's scaling relations~\cite{kesten} and with~\cite{lsw5}, has been used
by Smirnov and Werner \cite{sw} to obtain rigorously the values of several
critical exponents.


Consider a percolation model with distribution $P_p$ on a two-dimensional (regular)
lattice $\mathbb L$ such that the critical point $p_c$ is strictly between zero and
one. Let $C_0$ be the black cluster containing the origin and $|C_0|$ its cardinality,
then $\theta(p) = P_p(|C_0| = \infty)$ is the \emph{percolation probability}.
Arguments from theoretical physics suggest that, under general circumstances,
$\theta(p)$ behaves roughly like $(p-1/2)^{5/36}$ as $p$ approaches $p_c$ from above.
It is also believed that the \emph{connectivity function}
\begin{equation} \label{connect-function}
\tau_p(x) = P_p(\text{the origin and } x \text{ belong to the same cluster})
\end{equation}
behaves, for $0<p<p_c$ and $|x|$ large, 
like $\exp{(-|x|/\xi(p))}$, for some $\xi(p)$ satisfying
$\xi(p) \to \infty$ as $p \uparrow p_c$.
The \emph{correlation length} $\xi(p)$ is defined by
\begin{equation} \label{correlation-length}
\xi(p)^{-1} = \lim_{|x| \to \infty} \left\{ - \frac{1}{|x|} \log \tau_p(x) \right\}
\end{equation}
and is supposed to behave like $(p_c-p)^{-4/3}$ as $p \uparrow p_c$. The
\emph{mean cluster size} $\chi(p) = E_p(|C_0|)$ is also believed to diverge
with a power law behavior, $(p_c-p)^{-43/18}$, as $p \uparrow p_c$.


In~\cite{sw}, Smirnov and Werner proved the following result for the percolation
model on the triangular lattice discussed in this paper.
\begin{theorem} \label{thm-crit-exp}
With the above notation, we have
\begin{eqnarray}
\lim_{p \downarrow 1/2} \frac{\log \theta(p)}{\log (p - 1/2)} = 5/36, \\
\lim_{p \uparrow 1/2} \frac{\log \xi(p)}{\log (1/2 - p)} = - 4/3, \label{eq-corr-length} \\
\lim_{p \uparrow 1/2} \frac{\log \chi(p)}{\log (1/2 - p)} = - 43/18.
\end{eqnarray}
\end{theorem}

It remains a challenge to strengthen the above results and prove asymptotic
behavior up to constants. However, arguably the most important open problem is
obtaining (the existence of) critical exponents for percolation on other lattices.

\subsection{The Full Scaling Limit} \label{sec-full}

A single SLE$_6$ curve contains only limited information concerning connectivity
properties, and does not give a full description of the scaling limit. A more
complete description can be obtained in terms of loops, corresponding to the
scaling limit of cluster boundaries. Such loops should also be random and have
a conformally invariant distribution, closely related to SLE$_6$. This observation
motivated the work of Camia and Newman~\cite{cn1,cn2,cn4} and led
Werner~\cite{werner5,werner3} and Sheffield~\cite{sheffield} to define
and study the Conformal Loop Ensembles (CLEs) introduced in Sect.~\ref{sec-sle-cle}.

In~\cite{cn4}, it is shown that the collection of all cluster boundaries
contained in a Jordan domain $D$ converges in an appropriate sense to a
conformally invariant limit, called \emph{full scaling limit} in $D$. The
result can be stated as follows.

\begin{theorem} \label{main-result1}
In the continuum scaling limit, the probability distribution of the
collection of all boundary contours of critical site percolation on
the triangular lattice in a Jordan domain $D$ with monochromatic
boundary conditions converges to a probability distribution on
collections of continuous nonsimple loops in $D$. The limit is
conformally invariant in the following sense.
Let $D,D'$ be two Jordan domains and let $f:\overline D \to {\overline D}'$
a continuous function that maps $D$ conformally onto $D'$. Then the full
scaling limit in $D'$ is distributed like the image under $f$ of the full
scaling limit in $D$.
\end{theorem}

There is an explicit relation between the percolation full scaling
limit and SLE$_6$. It is shown in~\cite{cn4} that a process of loops
with the same distribution as in Theorem~\ref{main-result1} can be
constructed by a procedure in which each loop is obtained as the
concatenation of an SLE$_6$ path with (a portion of) another SLE$_6$
path. 
Moreover, as shown in the next theorem from~\cite{cn4}, the outermost
loops of the full scaling limit in a Jordan domain satisfy a conformal
restriction/renewal property, as in the definitions of the Conformal
Loop Ensembles of Werner~\cite{werner5} and Sheffield~\cite{sheffield}.
\begin{theorem} \label{thm-cle}
Consider the full scaling limit inside a Jordan domain $D$ and denote by ${\cal L}_D$
the collection of loops in $\overline D$ that are not surrounded by any other loop.
Consider an arc $\Gamma$ of $\partial D$ and let ${\cal L}_{D,\Gamma}$ be the set of
loops of ${\cal L}_D$ that touch $\Gamma$. Then, conditioned on ${\cal L}_{D,\Gamma}$,
for any connected component $D'$ of $D \setminus \overline{\cup \{ L:L \in {\cal L}_{D,\Gamma} \}}$,
the loops in $\overline{D'}$ form a random collection of loops distributed as an
independent copy of ${\cal L}_D$ conformally mapped to $D'$.
\end{theorem}

The full scaling limit in the whole plane can be obtained by taking
a sequence of domains $D$ tending to ${\mathbb R}^2$, as explained
in~\cite{cn1,cn2}. At the critical point, with probability one there
is no infinite black or white cluster, therefore the cluster boundaries
form finite (but arbitrarily large) loops even in the whole plane. The
following theorem from~\cite{cn2} lists some of the properties of the
full scaling limit in the plane.

\begin{theorem} \label{features}
The full scaling limit in the whole plane has the following
properties, which are valid with probability one:
\begin{enumerate}
\item It is a random collection of countably many
noncrossing continuous loops in the plane.
The loops can and do touch themselves and each other
many times, but there are no triple points; i.e. no three
or more loops can come together at the same point, and a
single loop cannot touch the same point more than twice,
nor can a loop touch a point where another loop touches
itself.
\item Any deterministic point $z$ in the plane (i.e.,
chosen independently of the loop process) is surrounded by
an infinite family of nested loops with diameters going to
both zero and infinity; any annulus about that point with
inner radius $r_1 > 0$ and outer radius $r_2 < \infty$
contains only a finite number $N(z, r_1, r_2)$ of those loops.
Consequently, any two distinct deterministic points of the
plane are separated by loops winding around each of them.
\item Any two loops are connected by a finite ``path''
of touching loops.
\end{enumerate}
\end{theorem}

\section{Near-Critical Percolation} \label{sec-near-crit}

Using the full scaling limit, one can attempt to understand the geometry of
near-critical scaling limits, where the percolation density $p$ tends to the
critical one $p_c$ in an appropriate way as the lattice spacing tends to zero:
\begin{equation} \label{eq-near-critical}
p=p_c+\lambda\delta^{\alpha}
\end{equation}
where $\delta$ is the lattice spacing, $\lambda \in (-\infty,\infty)$,
and $\alpha$ is set equal to $3/4$ to get nontrivial $\lambda$-dependence
in the limit $\delta \to 0$ (see, e.g., \cite{aizenman1,aizenman2,ab,bcks}).

A heuristic analysis~\cite{cfn1,cfn2} based on a natural ansatz leads
to a one-parameter family of loop models (i.e., probability measures
on random collections of loops), with the critical full scaling limit
corresponding to the special choice $\lambda=0$. Except for the latter
case, these measures are not scale invariant (and therefore are not
conformally invariant), but are \emph{conformally covariant} in an
appropriate sense and are mapped into one another by scale transformations.

The approach proposed in~\cite{cfn1,cfn2} is based on a ``Poissonian
marking" of double points of the critical full scaling limit, i.e.,
points where two loops of the critical full scaling limit touch each
other or a loop touches itself. These points correspond, back on the
lattice, to ``macroscopically pivotal" hexagons where four, long
${\cal T}$-paths of alternating colors come together.
The approach leads to a conceptual framework that can be used to define a
renormalization group flow (under the action of dilations), and to describe
the scaling limit of related models, such as invasion and dynamical percolation
and the minimal spanning tree. In particular, it helps explain why the scaling
limit of the minimal spanning tree may be scale invariant but \emph{not}
conformally invariant, as first observed numerically by Wilson~\cite{wi}.


Partly following some of the ideas and heuristic arguments put forward in~\cite{cfn1,cfn2},
significant progress has recently been made by Garban, Pete and Schramm in the rigorous
understanding of the near-critical scaling limit of percolation, in particular concerning
the uniqueness of the scaling limit and some of its properties, including rotation invariance
and conformal covariance.

Other interesting results are contained in~\cite{nowe}. There, the near-critical
scaling regime is not identified with the $\delta \to 0$ limit of a percolation
model with density of open sites given by~(\ref{eq-near-critical}) with $\alpha=3/4$,
but rather through the property that the correlation length introduced in
Sect.~\ref{sec-crit-exp} remains bounded away from $0$ and $\infty$ as $\delta \to 0$.
When that happens, it is shown in~\cite{nowe} that any subsequential scaling limit
of a percolation interface (between two points on the boundary of a bounded domain)
is mutually singular with respect to SLE$_6$, a somewhat surprising result, given
that other properties, like the Hausdorff dimension, are shared by critical and
near-critical interfaces in the scaling limit. It is also shown in~\cite{nowe} that
away from the near-critical regime described above the scaling limit of a percolation
interface is either SLE$_6$ (if the correlation length diverges) or degenerate (if
the correlation length goes to zero --- see~\cite{nowe} for more details).

In~\cite{cjm} it is shown that similar results for the collection of all boundaries
in the whole plane can be proved using ideas and tools that originate in Kesten's
seminal paper~\cite{kesten} (see also~\cite{nolin}). In particular, there are only
three possible types of scaling limits for the collection of all percolation interfaces:
the degenerate one consisting of no curves at all, the critical one corresponding to
the scaling limit of critical percolation, and one in which any deterministic point
in the plane is surrounded with probability one by a largest loop and by a countably
infinite family of nested loops with radii going to zero. All three cases occur.

The first case corresponds to the scaling limit of cluster boundaries in the
subcritical and supercritical phases. The second one arises when $p \to 1/2$ so
fast, as $\delta \to 0$, that the critical scaling limit is obtained. The last case
is particularly interesting because the scaling limit is nontrivial, like the critical
one, but unlike the critical one it is not scale invariant. This situation is described,
depending on the context, as {\em near-critical}, {\em off-critical} or {\em massive}
scaling limit (where ``massive" refers to the persistence of a macroscopic correlation
length, which should give rise to what is known in the physics literature as a
{\em massive field theory}). Near-critical scaling limits differ qualitatively from
the critical one at large scales, since in the latter there is no largest loop around
any point (see {\em 2.} in Theorem~\ref{features}). At the same time, they resemble
the critical scaling limit at short scales because of the presence of infinitely many
nested loops with radii going to zero around any given point.



\end{document}